\documentclass{amsart}
\usepackage[utf8]{inputenc}
\usepackage{lmodern}

\usepackage{extdash} 
\usepackage{amsmath,amsthm,amsfonts,amssymb, microtype}
\usepackage{xr-hyper}
\usepackage[colorlinks=true]{hyperref}
\usepackage{verbatim}
\usepackage{pdfsync} 
\usepackage{stmaryrd}
\usepackage{marvosym}
\usepackage{colonequals}
\usepackage{mathrsfs}
\usepackage{bbold}

\usepackage{mathtools}
\usepackage{enumitem}
\usepackage{tikz}
\usetikzlibrary{arrows,calc,positioning}

\usetikzlibrary{backgrounds}

\usepackage{extdash}
\usepackage{amsthm,amsfonts,amssymb,amsmath, microtype}
\usepackage{xr-hyper}
\usepackage[colorlinks=true]{hyperref}
\usepackage{verbatim}
\usepackage{pdfsync}
\usepackage{comment} 

\usetikzlibrary{circuits.logic.IEC}
\usetikzlibrary{cd}
\usepackage{xcolor}
\usetikzlibrary{shapes}
\usetikzlibrary{decorations.markings}

\tikzstyle heightone=[scale=.7,shift={(0,-.3)}]
\tikzstyle heightones=[scale=.8,xscale=.35,shift={(0,.1)}]
\tikzstyle heightoneonehalf=[scale=.9,shift={(0,-.2)}]
\tikzstyle heighttwo=[scale=.9,shift={(0,-.4)}]
\tikzstyle heighttwos=[scale=.5,xscale=.6,shift={(0,-.1)}]
\tikzstyle heightthree=[scale=.6,shift={(0,-.9)}]
\tikzstyle heightthrees=[scale=.4,xscale=.7,shift={(0,-.2)}]

\tikzstyle arrowstyle=[blue,semitransparent,scale=2]

\tikzstyle basiclabel=[draw=none,fill=none,shape=rectangle,inner sep=2pt,scale=.8]
\tikzstyle leftlabel=[basiclabel,anchor=east]
\tikzstyle rightlabel=[basiclabel,anchor=west]
\tikzstyle bottomlabel=[basiclabel,anchor=north]
\tikzstyle toplabel=[basiclabel,anchor=south]

\tikzstyle vertex=[circle,draw,fill=black,inner sep=1pt]
\tikzstyle ciliation=[circle,draw=none,fill=red,inner sep=1pt,semitransparent]
\tikzstyle ciliatednode=[vertex,pin={[pin distance=1mm,pin edge={semitransparent,red},ciliation]#1:{}}]

\tikzstyle vector=[black,thick,rectangle,draw=gray!50!yellow,top color=yellow!30,bottom color=black!10,scale=.8,inner sep=2pt]
\tikzstyle small vector=[vector,scale=.8]
\tikzstyle plain vector=[rectangle,draw=none,fill=white,scale=.7]

\tikzstyle my signal=[black,thick,signal,signal pointer angle=120,draw=blue!50,top color=blue!20,bottom color=black!10,scale=.8,inner sep=2pt]
\tikzstyle matrix=[my signal,signal from=south,signal to=north]
\tikzstyle reverse matrix=[my signal,signal from=north,signal to=south]

\tikzstyle small matrix=[matrix,scale=.7]
\tikzstyle reverse small matrix=[reverse matrix,scale=.7]
\tikzstyle matrix on edge=[small matrix,sloped,rotate=-90]
\tikzstyle reverse matrix on edge=[small matrix,sloped,rotate=90]

\tikzstyle trivalent=[very thick]
\tikzstyle dotdotdot=[decorate,decoration={markings,
    mark=at position .3 with{\node{.};},
    mark=at position .5 with {\node{.};},
    mark=at position .7 with {\node{.};}}]

\tikzstyle wavyup=[out=90,in=-90]
\tikzstyle wavydown=[out=-90,in=90]

\tikzstyle symmetrizer=[rectangle,fill=gray!10,draw=black]
\tikzstyle permutation=[symmetrizer]
\tikzstyle antisymmetrizer=[rectangle,fill=black,draw=black]
\tikzstyle symlabel=[draw=none,fill=none,black,scale=.8]
\tikzstyle asymlabel=[draw=none,fill=none,white,scale=.8]

\newcommand{\K}{\mathrm{K}}

\newcommand{\C}{\mathcal{C}}

\newcommand{\WW}{\mathrm{W}}
\newcommand{\CC}{\mathrm{C}}

\newcommand{\Mor}{\mathrm{Mor}}

\newcommand{\nocontentsline}[3]{}
\newcommand{\tocless}[2]{\bgroup\let\addcontentsline=\nocontentsline#1{#2}\egroup}

\newtheorem{theorem}{Theorem}[section]
\newtheorem{proposition}[theorem]{Proposition}
\newtheorem{lemma}[theorem]{Lemma}

\newtheorem*{theorem*}{Theorem}

\theoremstyle{definition}
\newtheorem{definition}[theorem]{Definition}
\newtheorem{lemma/definition}[theorem]{Definition/Lemma}
\newtheorem{example}[theorem]{Example}

\newtheorem{remark}[theorem]{Remark}

\begin{document}

\title{Algebraic $K_0$ for unpointed homotopy Categories}
\author{Felix K\"ung}
\begin{abstract}
We introduce the notion of Grothendieck heaps for unpointed Waldhausen categories and unpointed stable $\infty$-categories. This allows an extension of the studies of $\K_0$ to the homotopy category of unpointed topological spaces.

\end{abstract} 
\maketitle
\let\thefootnote\relax\footnotetext{The author is a Postdoctoral Research Fellow at the Universit\'e Libre de Bruxelles. The research in this work has been supported by the MIS (MIS/BEJ - F.4545.21) Grant from the FNRS and an ACR Grant from the Universit\'e Libre de Bruxelles.}
\section*{Introduction}

The Grothendieck group of an abelian category is one of the most studied invariants in mathematics \cite{Milnor1971,Brown1973,Quillen1975,Rosenberg1994}. The main idea of its construction is to interpret the objects as generators of a group and induce relations by short exact sequences. This idea has been generalized to the setting of stable $\infty$-categories \cite{Barwick2016} and Quillen exact categories \cite{Quillen1973} which then has been generalized to Waldhausen categories \cite{Waldhausen1985} by only requiring cofibrations and weak equivalences. There the idea is to replace the relations induced by short exact sequences by relations arising from exact and coexact triangles, exact sequences or cofiber sequences. These incarnations are very useful for categorification, as they allow the decategorification of derived categories and other types of homotopy categories \cite{Khovanov2016,Khovanov2009,Stroppel2005,Elias2016}.

As a heap is a natural generalization of groups, we defined in \cite{Kueng2023} the Grothendieck heap of a possibly unpointed category, hereby pushouts with a monomorphic leg take over the role of short exact sequences. Following this definition it was requested to develop the possible generalization of this Grothendieck heap to the setting of unpointed Waldhausen categories, respectively unpointed stable $\infty$-categories. In this paper we realise this request.

We start by recalling the basic constructions of the Grothendieck group for Waldhausen categories and stable $\infty$-categories. We then define the concept of unpointed Waldhausen categories and construct its Grothendieck heap. In Theorem~\ref{Theorem Heapy K0 recovers classical K0 for unpointed Waldhausen} we prove that in the case of a Waldhausen category our construction recovers the original definition. In the last section we repeat the same process for unpointed stable $\infty$-categories and prove in Theorem~\ref{Theorem Heapy K0 recovers classical K0 for unpointed stable} that our construction recovers the original definition if the category is pointed.

\section*{Acknowledgements}

We would like to thank G. Janssens and the referee of our previous paper for the encouragement to develop rigorously the ideas of the current paper. Furthermore I am grateful to the SL-Math institute for the awesome facilities that allowed this work to be completed.

\section{Preliminaries}

In this section we collect some basic notions about heaps and the classical algebraic $\K_0$-groups for Waldhausen categories respectively stable $\infty$-categories. 

\subsection{Heaps}
For convenience of the reader we recall the necessary basic definitions and Lemma. For more rigorous treatment of the basic notions on Heaps we refer to the work of T. Brzeziński \cite{Brzezinski2020}.

As a general rule of fist, one can think of a heap as a group without a zero object.

\begin{definition}
A heap is a set $H$ together with a ternary operation 
\begin{align*}
H\times H\times H &\to H \\
\left(x,y,z\right)&\mapsto \left[x,y,z\right]
\end{align*} such that the following two relations hold:
\begin{align*}
\left[x,x,y\right]=&y=\left[y,x,x\right] &\text{unitality}\\
\left[x,y,\left[x',y',z'\right]\right]&=\left[\left[x,y,x'\right],y',z'\right] &\text{associativity.}
\end{align*}
\end{definition}

\begin{example}
Let $G$ be a group, then $G$ admits a heap structure given by $$\left[x,y,z\right]:=xy^{-1}z.$$
\end{example}

We can do an inverse construction and build a group out of a heap by forcing a zero object.

\begin{definition}
Let $H$ be a heap and $e\in H$, then the retract of $H$ along $e$ is the group $\left(H,+_e\right)$ given by the set $H$ and group structure $$x+_e y:= \left(x,e,y\right)$$
\end{definition}

\begin{definition}
A map $H\to H'$ is a morphism of heaps, if it is compatible with the heap bracket i.e.
$$f\left[x,y,z\right]=\left[fx,fy,fz\right].$$
\end{definition}

\begin{lemma}\cite[Lemma~1.18]{Kueng2023}\label{reconstruct heap from retract}
Let $H$ be a heap and $e \in H$. Then we have 
$$\left[a,b,c\right]=a-_e b +_e c.$$
\end{lemma}

\subsection{Algebraic $\K_0$ for Waldhausen categories}

We recall the following definitions and constructions in order to later adapt them on the non-pointed level. The notions we cover here are very classical and there are countless references for them, for instance \cite{Weibel2013,Waldhausen1985}.

\begin{definition}
A Waldhausen category is a triple $\left(\C,\CC,\WW\right)$ consisting of a category $\C$, a subcategory called cofibrations $\CC\subset \C$ and a class of morphisms called weak equivalences $\WW\subset \Mor\left(\C\right)$. Such that the following axioms are satisfied:
\begin{itemize}
\item Isomorphisms are cofibrations.
\item There is a zero object $0\in \C$ and $0 \rightarrowtail X$ is a cofibration for all $X\in \C$.
\item If $Y\rightarrowtail Z$ is a cofibration and $Y\to X$ any morphism, then the pushout $X\sqcup_Y Z$ along these maps exists and $X\rightarrowtail X\sqcup_Y Z$ is a cofibration.
\item Isomorphisms are weak equivalences.
\item Weak equivalences are closed under composition.
\item Given a commutative diagram of the shape 

$$\tikz[heighttwo,xscale=1,yscale=1,baseline]{
\node (Y) at (0,4) {$Y$};
\node (X) at (0,2) {$X$};
\node (Z) at (2,4) {$Z$};
\node (Y1) at (2,2) {$Y'$};
\node (X1) at (2,0) {$X'$};
\node (Z1) at (4,2) {$Y'$};

\draw[->]
(X) edge node[below left] {$h$}(X1)
(Y)edge node[above right] {$w$}(Y1)
(Z)edge node[above right] {$s$}(Z1)
(Y)edge node[left] {$f$}(X)
(Y1) edge node[right] {$f'$}(X1)
;
\draw[>->]
(Y)edge node[above] {$c$} (Z)
(Y1)edge node[below] {$c'$} (Z1)
;
}
$$
with $c,c' \in \CC$ and $h,w,s \in \WW$, then the induced morphism of pushouts $X\sqcup_Y Z \xrightarrow{\sim} X'\sqcup_{Y'} Z'$ is a weak equivalence.
\end{itemize}
\end{definition}

We usually will write $X\rightarrowtail Y$ for cofibrations, $X \xrightarrow{\sim} Y$ for weak equivalences and write $\C$ for a Waldhausen category instead of the full triple $\left(\C,\CC,\WW\right)$ for improved readability.

\begin{remark}
The axiom of $0\rightarrowtail X$ being a cofibration is used in order to make the fundamental group of the Waldhausen $S$-construction see all objects and hence can be thought of as a connectivity axiom. As we want to work with unpointed categories, it is natural to not require the spaces to be connected.
\end{remark}

\begin{definition}
A sequence $X \rightarrowtail Y \to X/Y$ is a cofiber sequence if there is a pushout diagram of the shape

$$\tikz[heighttwo,xscale=1,yscale=1,baseline]{
\node (Y) at (0,2) {$Y$};
\node (X) at (0,0) {$0$};
\node (Z) at (2,2) {$X$};
\node (P) at (2,0) {$X/Y$.};
\node (po)at (1,1) {$\ulcorner$};
\draw[->]
(Y) edge (X)
(Z)edge (P)
;
\draw[>->]
(Y)edge (Z)
(X) edge (P)
;
}
$$
\end{definition}

\begin{remark}
In the following Definition~\ref{Definition classical K0 for Waldhausen} we require $W$ to have a set of equivalence classes in order for the generating elements to define a set.
\end{remark}

\begin{definition}\label{Definition classical K0 for Waldhausen}
Let $\left(\C,\CC,\WW \right)$ be a Waldhausen category with a set of weak equivalence classes. Then the Grothendieck group $\K_0\left(\C\right)$ is defined to be the free abelian group generated by the weak equivalence classes $\overline{X}$ and relations induced by cofiber sequences $X\rightarrowtail Y \to Y/X$ i.e.
$$\K_0\left(\C\right):=\left\langle \overline{X} \middle| \overline{X} + \overline{X/Y} = \overline{Y} \text{ for a cofiber sequence }X\rightarrowtail Y \to Y/X\right\rangle .$$
\end{definition}

\begin{definition}
A functor $F: \C \to \C'$ between two Waldhausen categories $\left(\C,\CC,\WW\right),\left(\C',\CC',\WW'\right)$ is exact if

\begin{itemize}
\item it is pointed i.e. $F \left(0\right)\cong 0$,
\item preserves cofibrations i.e. $F \left(c\right) \in \CC'$ if $c\in \CC$,
\item preserves weak equivalences i.e. $F\left(w\right) \in \WW'$ if $w\in W$
\item and preserves pushouts along cofibrations i.e. the canonical morphism $F\left(X\right) \sqcup_{F\left(Y\right)} F\left(Z\right) \to F\left(X\sqcup_Y Z\right)$ is an isomorphism for all cofibrations $Y\rightarrowtail Z$.
\end{itemize}
\end{definition}

\begin{proposition}
Let $F: \C \to \C' $ be an exact functor of Waldhausen categories. Then $F$ induces a morphism of groups
\begin{align*}
F:\K_0\left(\C\right) &\to \K_0\left(\C'\right)\\
\overline{X}&\mapsto \overline{F \left(X\right)}.
\end{align*}
\end{proposition}

\subsection{Algebraic $\K_0$ for stable $\infty$-categories}

In the following we discuss $\infty$-categories, for reference we refer the reader to \cite{Lurie2006,Lurie2009}. However, as we only consider $\K_0$ here, knowing all the heavy machinery is not essential.

\begin{definition}
Let $\C$ be a pointed $\infty$-category, then we recall the following notions.
\begin{itemize}
\item A cokernel of $f: X\to Y$ in $\C$ is a pushout of the form
$$\tikz[heighttwo,xscale=1,yscale=1,baseline]{
\node (Y) at (0,2) {$X$};
\node (X) at (0,0) {$0$};
\node (Z) at (2,2) {$Y$};
\node (P) at (2,0) {$Z$.};
\node (po) at (1,1) {$\ulcorner$};

\draw[->]
(Y) edge (X)
(Z)edge (P)
(Y)edge node[above] {$f$} (Z)
(X) edge (P)
;
}
$$
\item A kernel of $f: P\to Z$ in $\C$ is a pullback of the form
$$\tikz[heighttwo,xscale=1,yscale=1,baseline]{
\node (Y) at (0,2) {$X$};
\node (X) at (0,0) {$0$};
\node (Z) at (2,2) {$Y$};
\node (P) at (2,0) {$Z$.};
\node (po) at (1,1) {$\lrcorner$};

\draw[->]
(Y) edge (X)
(Z)edge node[right] {$f$}(P)
(Y)edge (Z)
(X) edge (P)
;
}
$$
\item An exact triangle $X\to Y\to Z$ in $\C$ is a pushout of the form 
$$\tikz[heighttwo,xscale=1,yscale=1,baseline]{
\node (Y) at (0,2) {$X$};
\node (X) at (0,0) {$0$};
\node (Z) at (2,2) {$Y$};
\node (P) at (2,0) {$Z$.};
\node (po) at (1,1) {$\ulcorner$};

\draw[->]
(Y) edge (X)
(Z)edge (P)
(Y)edge (Z)
(X) edge (P)
;
}
$$
\item A coexact triangle $X\to Y\to Z$ in $\C$ is a pullback of the form
$$\tikz[heighttwo,xscale=1,yscale=1,baseline]{
\node (Y) at (0,2) {$X$};
\node (X) at (0,0) {$0$};
\node (Z) at (2,2) {$Y$};
\node (P) at (2,0) {$Z$.};
\node (po) at (1,1) {$\lrcorner$};

\draw[->]
(Y) edge (X)
(Z)edge (P)
(Y)edge (Z)
(X) edge (P)
;
}
$$
\end{itemize}
\end{definition}

\begin{definition}\cite[Definition~2.9]{Lurie2006}
A stable $\infty$-category $\C$ is an $\infty$-category $\C$ such that
\begin{itemize}
\item $\C$ is pointed,
\item every morphism $f: X \to Y$ admits a kernel and cokernel,
\item a triangle is exact if and only if its coexact.
\end{itemize}
\end{definition}

However, we will use the following equivalent characterization which allows a straightforward generalization to the unpointed setting.

\begin{lemma}\cite[Proposition~4.4]{Lurie2006}\label{Lemma stable infinity categories are pointed admitting bicartesian squares}
An $\infty$-category $\C$ is stable if and only if
\begin{itemize}
\item $\C$ is pointed,
\item $\C$ admits finite limits and colimits,
\item a diagram in $\C$
$$\tikz[heighttwo,xscale=1,yscale=1,baseline]{
\node (Y) at (0,2) {$Y$};
\node (X) at (0,0) {$X$};
\node (Z) at (2,2) {$Y'$};
\node (P) at (2,0) {$X'$};
\node (po) at (1,1) {$\square$};
\draw[->]
(Y) edge (X)
(Z)edge (P)
(Y)edge (Z)
(X) edge (P)
;
}
$$
is a pushout if and only if it is a pullback.
\end{itemize}
\end{lemma}

\begin{definition}
A diagram of the shape 
$$\tikz[heighttwo,xscale=1,yscale=1,baseline]{
\node (Y) at (0,2) {$Y$};
\node (X) at (0,0) {$X$};
\node (Z) at (2,2) {$Y'$};
\node (P) at (2,0) {$X'$.};
\node (po) at (1,1) {$\square$};
\draw[->]
(Y) edge (X)
(Z)edge (P)
(Y)edge (Z)
(X) edge (P)
;
}
$$
that is both a pushout and a pullback is called bicartesian.
\end{definition}

\begin{definition}
Let $\C$ be a stable $\infty$-category then its Grothendieck group $\K_0\left(\C\right)$ is the free abelian group of equivalence classes of objects with relations induced by exact triangles:
$$\K_0\left(\C\right):=\left\langle \overline{X}\middle| \overline{X} +\overline{Z}=\overline{Y} \text{ if there is an exact triangle }X\to Y\to Z\right\rangle .$$
\end{definition}

\section{$\K_0$ for unpointed homotopy categories admitting coproducts}

We will start by considering unpointed Waldhausen categories and then continue to unpointed stable $\infty$-categories. In both cases we need to be careful to not invite an Eilenberg swindle. However, both cases have their own way of elegantly avoiding that.

\subsection{$\K_0$ for unpointed Waldhausen categories}

We start by defining the unpointed incarnation of Waldhausen categories. Observe that the property of the unique morphism from the zero object being a cofibration is mainly used to make sure that the Waldhausen $S$/$T$-construction indeed reaches all objects. In particular, as we only care about the naive $\K_0$ we can ignore this issue.

\begin{definition}
An unpointed Waldhausen category is a triple $\left(\C,\CC,\WW\right)$ consisting of a category $\C$, a subcategory of cofibrations $\CC\subset \C$ and a class of morphisms called weak equivalences $\WW\subset \Mor\left(\C\right)$. Such that the following axioms are satisfied:
\begin{itemize}
\item Isomorphisms are cofibrations.
\item If $Y\rightarrowtail Z$ is a cofibration and $Y\to X$ any morphism, then the pushout $X\sqcup_Y Z$ along these maps exists and $X\rightarrowtail X\sqcup_Y Z$ is a cofibration.
\item Isomorphisms are weak equivalences.
\item Weak equivalences are closed under composition.
\item Given a commutative diagram of the shape 

$$\tikz[heighttwo,xscale=1,yscale=1,baseline]{
\node (Y) at (0,4) {$Y$};
\node (X) at (0,2) {$X$};
\node (Z) at (2,4) {$Z$};
\node (Y1) at (2,2) {$Y'$};
\node (X1) at (2,0) {$X'$};
\node (Z1) at (4,2) {$Y'$};

\draw[->]
(X) edge node[below left] {$h$}(X1)
(Y)edge node[above right] {$w$}(Y1)
(Z)edge node[above right] {$s$}(Z1)
(Y)edge node[left] {$f$}(X)
(Y1) edge node[right] {$f'$}(X1)
;
\draw[>->]
(Y)edge node[above] {$c$} (Z)
(Y1)edge node[below] {$c'$} (Z1)
;
}
$$
with $c,c' \in \CC$ and $h,w,s \in \WW$, then the induced morphism of pushouts $X\sqcup_Y Z \xrightarrow{\sim} X'\sqcup_{Y'} Z'$ is a weak equivalence.
\end{itemize}
\end{definition}

\begin{remark}
In the following we, similarly to the classical case, assume that $\WW$ has a set of weak equivalence classes in order to guarantee that the generating objects form a set.
\end{remark}

\begin{definition}
Let $\left(\C,\CC,\WW \right)$ be an unpointed Waldhausen category with a set of weak equivalence classes. Then the Grothendieck heap $\K_0\left(\C\right)$ is defined to be the free abelian heap generated by the weak equivalence classes $\overline{X}^\WW$ and relations induced by pushouts i.e.
$$\K_0\left(\C\right):=\left\langle \overline{X} \middle| \left[\overline{X} , \overline{Y},\overline{Z}\right] = \overline{X\sqcup_Y Z}\text{ for } Y\rightarrowtail Z \text{ or } Y\rightarrowtail X \right\rangle .$$
\end{definition}

\begin{remark}
Observe that similar to the classical cases and the cases treated in \cite{Kueng2023} the definition of a Waldhausen category only considering pushouts along cofibrations gives a form of minimality and hence prevents trivial Eilenberg swindle.
\end{remark}

\begin{theorem}\label{Theorem Heapy K0 recovers classical K0 for unpointed Waldhausen}
Let $\left(\C,\CC,\WW\right)$ be a Waldhausen category. Then we have that $\K_0 \left(\C\right)$ recovers the original definition of $\K_0$ by considering the retract along the zero object $\left(\K_0\left(\C\right),+_{\overline{0}}\right)$ i.e. 
$$\overline{X} +_{\overline{0}} \overline{Y} =\left[ \overline{A} , \overline{0} , \overline{B} \right].$$
\begin{proof}
Throughout this proof we will distinguish the a priori different constructions by marking the classical constructions with a subscript $c$.

Consider the map 
\begin{align*}
\psi:\K_0\left(\C\right)&\to \K_0\left(\C\right)_c\\
\overline{X}&\mapsto \overline{X}_c.
\end{align*}

We need to first show that this is well defined. Consider a generating relation of $\K_0\left(\C\right)$: Let 
\begin{equation*}
\tikz[heighttwo,xscale=2,yscale=2,baseline]{
\node (X1) at (0,1) {$Y$};
\node (X2) at (0,0) {$X$};
\node (Y1) at (1,1) {$Z$};
\node (X+Y) at (1,0) {$W$,};
\node (po)at (0.5,0.5) {$\ulcorner$};

\draw[>->]
(X1) edge  (X2);
\draw[->]
(X1) edge  (Y1)
(X2) edge  (X+Y)
(Y1)edge (X+Y);
}
\end{equation*}
be a pushout inducing $\left[\overline{X},\overline{Y},\overline{Z}\right]=\overline{W}$. This becomes after passing to $\left(\K_0\left(\C\right),+_{\overline{0}}\right)$ via Lemma~\ref{reconstruct heap from retract}
\begin{align*}
\left[\overline{X},\overline{Y},\overline{Z}\right]&=\overline{W}\\
\overline{X}-_{\overline{0}} \overline{Y}+_{\overline{0}}\overline{Z}&=\overline{W}\\
\overline{X} +_{\overline{0}} \overline{Z} &=\overline{Y} +_{\overline{0}} \overline{W}.
\end{align*} 

Now consider the diagram

\begin{equation*}
\tikz[heighttwo,xscale=2,yscale=2,baseline]{
\node (X1) at (0,2) {$Y$};
\node (X2) at (0,1) {$X$};
\node (Y1) at (1,2) {$Z$};
\node (X+Y) at (1,1) {$W$};
\node (0) at (0,0) {$0$};
\node (X+Y1) at (1,0) {$W.$};
\node (po)at (0.5,0.5) {$\ulcorner$};
\node (po)at (0.5,1.5) {$\ulcorner$};

\draw[>->]
(X1) edge  (X2)
(X2) edge   (X+Y)
(0) edge   (X+Y1);
\draw[->]
(X1) edge  (Y1)
(X2) edge  (X+Y)
(X+Y) edge  (X+Y1)
(X2) edge  (0)
(Y1)edge (X+Y);
}
\end{equation*}

By pasting of pushouts the total diagram is a pushout too. In particular we have in $\K_0\left(\C\right)_c$ the relation
$$\overline{Y}_c+\overline{W'}_c=\overline{Z}_c.$$
Now as the bottom square is a pushout as well we have $\overline{W'}_c=\overline{W}_c-\overline{X}_c$, and substituting gives the equation
$$\psi\overline{Y} + \psi\overline{C}=\overline{Y}_c+\overline{W}_c=\overline{Z}_c+\overline{X}_c=\psi\overline{Z}+\psi\overline{W}.$$
So $\psi$ is well defined.

Furthermore $\psi$ defines a group morphism as we have, since $0\rightarrowtail X$ is a cofibration for all $X\in \C$.
\begin{align*}
\psi\left(\overline{X} +_{\overline{0}} \overline{Y}\right)&=\psi\left(\left[\overline{X},\overline{0},\overline{Y}\right]\right)&*+_{\overline{0}}* = [*,\overline{0},*]\\
&=\psi\left({\overline{X\sqcup_{0}Y}}\right) &[*,*,*]=*\sqcup_* *\\
&=\psi\left({\overline{X\oplus Y}}\right) &\text{pushout along 0 is } \oplus \\
&=\overline{X\oplus Y}_c &\psi\left(\overline{X}\right)=\overline{X}_c\\
&=\overline{X}_c+\overline{Y}_c,
\end{align*}
where the last equality holds by the cofibersequence $X\rightarrowtail X\oplus Y \to Y$, and $X\rightarrowtail X\oplus Y$ is induced by the pushout

$$\tikz[heighttwo,xscale=2,yscale=2,baseline]{
\node (X1) at (0,1) {$0$};
\node (X2) at (0,0) {$X$};
\node (Y1) at (1,1) {$ Y$};
\node (X+Y) at (1,0) {$X\oplus Y$.};
\node (po)at (0.5,0.5) {$\ulcorner$};

\draw[>->]
(X1) edge  (Y1)
(X2) edge  (X+Y);
\draw[->]
(X1) edge  (X2)
(Y1)edge (X+Y);
}$$

By construction $\psi$ is surjective, in particular we now need to show injectivity. 

To prove $\ker\left(\psi\right)=0$ it is enough to show that for a cofiber sequence $ X \rightarrowtail Y \to Z $ we have $\overline{Y}=\overline{X} +_{\overline{0}}\overline{Z}$. Note that this is equivalent to
\begin{align*}
 \overline{Y}&=\overline{X}  +_{\overline{0}} \overline{Z}\\
\overline{Y}-_{\overline{0}} \overline{X}&=\overline{Z}\\
 \overline{0} +_{\overline{0}}\overline{Y}-_{\overline{0}} \overline{X} &=\overline{Z} \iff
\left[\overline{0},\overline{Y},\overline{X}\right]=\overline{Z}.& \text{Lemma~\ref{reconstruct heap from retract}}
\end{align*}
But this holds as by definition a cofiber sequence is given by a pushout
$$\tikz[heighttwo,xscale=2,yscale=2,baseline]{
\node (X1) at (0,1) {$X$};
\node (X2) at (0,0) {$0$};
\node (Y1) at (1,1) {$Y$};
\node (X+Y) at (1,0) {$Z$,};
\node (po)at (0.5,0.5) {$\ulcorner$};

\draw[>->]
(X1)edge  (Y1)
(X2) edge (X+Y);
\draw[->]
(X1) edge  (X2)
(Y1)edge (X+Y);
}
$$
and so $\psi$ is injective and hence an isomorphism.
\end{proof}	
\end{theorem}

\begin{definition}
A functor $F: \C \to \C'$ between two unpointed Waldhausen categories $\left(\C,\CC,\WW\right),\left(\C',\CC',\WW'\right)$ is exact if it

\begin{itemize}
\item preserves cofibrations i.e. $F c \in \CC'$ if $c\in \CC$,
\item preserves weak equivalences i.e. $Fw \in \WW'$ if $w\in W$
\item and preserves pushouts along cofibrations i.e. the canonical morphism $FX \sqcup_{FY} FZ \to F\left(X\sqcup_Y Z\right)$ is an isomorphism for all cofibrations $Y\rightarrowtail Z$.
\end{itemize}
\end{definition}

\begin{proposition}
Let $F: \C \to \C' $ be an exact functor of unpointed Waldhausen categories. Then $F$ induces a morphism of heaps
\begin{align*}
F:\K_0\left(\C\right) &\to \K_0\left(\C'\right)\\
\overline{X}&\mapsto \overline{F \left(X\right)}.
\end{align*}
\begin{proof}
As the functor $F$ maps weak equivalences to weak equivalences it preserves weak equivalence classes of objects and we need to check that $\K_0\left(F\right)$ respects the generating relations. For this consider a pushout giving rise to a generating relation in $\K_0\left(\C\right)$:
$$\tikz[heighttwo,xscale=2,yscale=2,baseline]{
\node (X1) at (0,1) {$Y$};
\node (X2) at (1,1) {$X$};
\node (Y1) at (0,0) {$Z$};
\node (X+Y) at (1,0) {$W$.};
\node (po)at (0.5,0.5) {$\ulcorner$};

\draw[->]
(X1)edge  (Y1)
(X2) edge (X+Y)
(Y1) edge (X+Y);
\draw[>->]
(X1) edge  (X2);
}
$$
Then, since $F$ preserves cofibrations and pullbacks along cofibrations, we have the pushout 
$$\tikz[heighttwo,xscale=2,yscale=2,baseline]{
\node (X1) at (0,1) {$FY$};
\node (X2) at (1,1) {$FX$};
\node (Y1) at (0,0) {$FZ$};
\node (X+Y) at (1,0) {$FW$.};
\node (po)at (0.5,0.5) {$\ulcorner$};

\draw[->]
(X1)edge  (Y1)
(X2) edge (X+Y)
(Y1) edge (X+Y);
\draw[>->]
(X1) edge  (X2);
}
$$
This by definition means that $[\overline{FX},\overline{FY},\overline{FZ}]=\overline{FW}$ as claimed.
\end{proof}
\end{proposition}

\subsection{$\K_0$ for unpointed stable $\infty$-categories}

In the case of $\infty$-categories one can use an alternative approach for guaranteeing minimality of pushouts and preventing Eilenberg-swindle, bicartesian squares. 

 We use the characterization by Lemma~\ref{Lemma stable infinity categories are pointed admitting bicartesian squares} of stable $\infty$-categories as this allows straightforward generalization to the unpointed case.

Due to the depth of the theory of $\infty$-categories and ease of reading we restrict our exposition to the definition of the heap $\K_0$ for unpointed stable $\infty$-categories and that it indeed recovers the classical definition.

\begin{definition}
An unpointed stable $\infty$-category $\C$ is an $\infty$-category $\C$ such that
\begin{itemize}
\item $\C$ admits finite limits and colimits,
\item a diagram in $\C$
$$\tikz[heighttwo,xscale=1,yscale=1,baseline]{
\node (Y) at (0,2) {$Y$};
\node (X) at (0,0) {$X$};
\node (Z) at (2,2) {$Y'$};
\node (P) at (2,0) {$X'$};
\node (po) at (1,1) {$\square$};
\draw[->]
(Y) edge (X)
(Z)edge (P)
(Y)edge (Z)
(X) edge (P)
;
}
$$
is a pushout if and only if it is a pullback.
\end{itemize}
\end{definition}

\begin{definition}
The Grothendieck heap of an unpointed stable $\infty$-category is the free abelian heap of equivalence classes of objects with relations induced by bicartesian squares:
$$\K_0\left(\C\right):=\left\langle \overline{X}\middle| \left[\overline{X},\overline{Y},\overline{Z}\right]=\overline{X\sqcup_Y Z}\right\rangle .$$
\end{definition}

\begin{theorem}\label{Theorem Heapy K0 recovers classical K0 for unpointed stable}
Let $\C$ be a stable $\infty$-category. Then we have that $\K_0 \left(\C\right)$ recovers the original definition of $\K_0$ by considering the retract along the zero object $\left(\K_0\left(\C\right),+_{\overline{0}}\right)$ i.e. 
$$\overline{X} +_{\overline{0}} \overline{Y} =\left[ \overline{A} , \overline{0} , \overline{B} \right].$$
\begin{proof}
As above we distinguish the a priori different constructions by marking the classical constructions with a subscript $c$.

Consider the map
\begin{align*}
\psi:\K_0\left(\C\right)&\to \K_0\left(\C\right)_c\\
\overline{X}&\mapsto \overline{X}_c.
\end{align*}

We need to first show that $\psi$ is well defined. Consider for this a generating relation of $\K_0\left(\C\right)$: Let 
\begin{equation*}
\tikz[heighttwo,xscale=2,yscale=2,baseline]{
\node (X1) at (0,1) {$Y$};
\node (X2) at (0,0) {$X$};
\node (Y1) at (1,1) {$Z$};
\node (X+Y) at (1,0) {$W$};
\node (po)at (0.5,0.5) {$\square$};

\draw[->]
(X1) edge  (X2)
(X1) edge  (Y1)
(X2) edge  (X+Y)
(Y1)edge (X+Y);
}
\end{equation*}
be a bicartesian square inducing $\left[\overline{X},\overline{Y},\overline{Z}\right]=\overline{W}$. This becomes after passing to $\left(\K_0\left(\C\right),+_{\overline{0}}\right)$ via Lemma~\ref{reconstruct heap from retract}
\begin{align*}
\left[\overline{X},\overline{Y},\overline{Z}\right]&=\overline{W}\\
\overline{X}-_{\overline{0}} \overline{Y}+_{\overline{0}}\overline{Z}&=\overline{W}\\
\overline{X} +_{\overline{0}} \overline{Z} &=\overline{Y} +_{\overline{0}} \overline{W}.
\end{align*} 

Now consider the diagram

\begin{equation*}
\tikz[heighttwo,xscale=2,yscale=2,baseline]{
\node (X1) at (0,2) {$Y$};
\node (X2) at (0,1) {$X$};
\node (Y1) at (1,2) {$Z$};
\node (X+Y) at (1,1) {$W$};
\node (0) at (0,0) {$0$};
\node (X+Y1) at (1,0) {$W'$.};
\node (po)at (0.5,0.5) {$\square$};
\node (po)at (0.5,1.5) {$\square$};

\draw[->]
(X1) edge  (X2)
(X2) edge   (X+Y)
(0) edge   (X+Y1)
(X1) edge  (Y1)
(X2) edge  (X+Y)
(X+Y) edge  (X+Y1)
(X2) edge  (0)
(Y1)edge (X+Y);
}
\end{equation*}

By pasting of bicartesian squares the total diagram is bicartesian too, this gives rise to an exact triangle $Y\to Z \to W'$. In particular we have the relation in $\K_0\left(\C\right)_c$
$$\overline{Y}_c+\overline{W'}_c=\overline{Z}_c.$$
Now as the bottom square is bicartesian we have $\overline{W'}_c=\overline{W}_c-\overline{X}_c$, and substituting gives the equation
$$\psi\overline{Y} + \psi\overline{C}=\overline{Y}_c+\overline{W}_c=\overline{Z}_c+\overline{X}_c=\psi\overline{Z}+\psi\overline{W}.$$
So $\psi$ is well defined.

Furthermore $\psi$ defines a group morphism as we have:
\begin{align*}
\psi\left(\overline{X} +_{\overline{0}} \overline{Y}\right)&=\psi\left(\left[\overline{X},\overline{0},\overline{Y}\right]\right)&*+_{\overline{0}}* = [*,\overline{0},*]\\
&=\psi\left({\overline{X\sqcup_{0}Y}}\right) &[*,*,*]=*\sqcup_* *\\
&=\psi\left({\overline{X\oplus Y}}\right) &\text{pushout along 0 is } \oplus \\
&=\overline{X\oplus Y}_c &\psi\left(\overline{X}\right)=\overline{X}_c\\
&=\overline{X}_c+\overline{Y}_c,
\end{align*}
where the last equality holds by the exact triangle $X\rightarrowtail X\oplus Y \to Y$.

By construction $\psi$ is surjective, in particular it remains to show injectivity. 

To prove $\ker\left(\psi\right)=0$ it is enough to show that for an exact triangle $ X \to Y \to Z $ we have $\overline{Y}=\overline{X} +_{\overline{0}}\overline{Z}$. Note that this is equivalent to
\begin{align*}
 \overline{Y}&=\overline{X}  +_{\overline{0}} \overline{Z}\\
\overline{Y}-_{\overline{0}} \overline{X}&=\overline{Z}\\
 \overline{0} +_{\overline{0}}\overline{Y}-_{\overline{0}} \overline{X} &=\overline{Z} \iff
\left[\overline{0},\overline{Y},\overline{X}\right]=\overline{Z}.& \text{Lemma~\ref{reconstruct heap from retract}}
\end{align*}
But this holds as by definition an exact triangle is given by a bicartesian square
$$\tikz[heighttwo,xscale=2,yscale=2,baseline]{
\node (X1) at (0,1) {$X$};
\node (X2) at (0,0) {$0$};
\node (Y1) at (1,1) {$Y$};
\node (X+Y) at (1,0) {$Z$,};
\node (po)at (0.5,0.5) {$\square$};

\draw[->]
(X1)edge  (Y1)
(X2) edge (X+Y)
(X1) edge  (X2)
(Y1)edge (X+Y);
}
$$
and so $\psi$ is injective and hence an isomorphism.
\end{proof}	
\end{theorem}

\begin{remark}
A recent approach by J. Campbell, J. Kuijper, M. Merling and I. Zakharevich \cite{Campbell2023} studies the algebraic $\K$-theory of categories of squares. This approach together with an insight by Vokřínek \cite{Vokrinek2013a} that morphism spaces in unpointed stable model categories admit naturally the structure of a heap should give rise to an unpointed Waldhausen $T$-construction.
\end{remark}

 \bibliographystyle{hep}
\bibliography{algebraic_K0_for_unpointed_homotopy_Categories}  

\end{document}